# EXPLOSION PHENOMENA IN STOCHASTIC COAGULATION–FRAGMENTATION MODELS

By Wolfgang Wagner

*Weierstrass Institute for Applied Analysis and Stochastics*

First we establish explosion criteria for jump processes with an arbitrary locally compact separable metric state space. Then these results are applied to two stochastic coagulation–fragmentation models—the direct simulation model and the mass flow model. In the pure coagulation case, there is almost sure explosion in the mass flow model for arbitrary homogeneous coagulation kernels with exponent bigger than 1. In the case of pure multiple fragmentation with a continuous size space, explosion occurs in both models provided the total fragmentation rate grows sufficiently fast at zero. However, an example shows that the explosion properties of both models are not equivalent.

**1. Introduction.** Coagulation–fragmentation models are used in different application fields ranging from chemical engineering (reacting polymers, soot formation) or aerosol technology to astrophysics (formation of stars and planets). These models describe the behavior of a system of particles that are characterized by their sizes and move in a certain medium. The size of a particle changes either by coagulation (merging with another particle) or by fragmentation (breakage, splitting into pieces). We refer to the survey paper [1] for more details and references.

Deterministic coagulation–fragmentation models are nonlinear evolution equations governing the macroscopic behavior of the particle system. In the case of a discrete size variable and spatial homogeneity, the coagulation equation takes the form

$$\frac{\partial}{\partial t} c(t,x) = \frac{1}{2} \sum_{y=1}^{x-1} K(x-y,y) c(t,x-y) c(t,y)$$

(1.1)

$$-\sum_{y=1}^{\infty} K(x,y)c(t,x)c(t,y),$$

where $t \geq 0$ and $x = 1, 2, \ldots$. Equation (1.1) goes back to Smoluchowski [29] and describes the time evolution of the particle number density (relative number of particles of a given size). The coagulation kernel $K$ is determined by the physics of the driving medium. It was observed that a certain kind of phase transition occurs when $K$ grows sufficiently fast in its arguments. This phenomenon is called gelation and corresponds to a loss of mass in the macroscopic equation. The gelation point is defined as

(1.2) $\quad t_{\mathrm{gel}} = \inf\{t \geq 0 : m_1(t) < m_1(0)\} \qquad \text{where } m_1(t) = \sum_{x=1}^{\infty} x c(t,x).$

Gelation is interpreted as the formation of infinitely large clusters in finite time. The presence of such clusters is not reflected in (1.1).

Stochastic coagulation models go back to [14, 20, 22]. They are based on systems of particles

(1.3) $\qquad x_i^{(n)}(t), \qquad i = 1, \ldots, N^{(n)}(t), \qquad N^{(n)}(0) = n,$

which coagulate according to appropriate rates determined by the kernel $K$. These systems approximate the solution of the coagulation equation (1.1) in the sense

(1.4) $\qquad \sum_{x=1}^{\infty} \varphi(x) c(t,x) = \lim_{n \to \infty} \frac{1}{n} \sum_{i=1}^{N^{(n)}(t)} \varphi(x_i^{(n)}(t)), \qquad t \geq 0,$

for appropriate test functions $\varphi$ (sequences in the discrete case). We refer to [8, 16, 26, 27] concerning rigorous results. An alternative stochastic coagulation model is related to the mass flow equation

(1.5)
$$\frac{\partial}{\partial t} \tilde{c}(t,x) = \sum_{y=1}^{x-1} \frac{K(x-y,y)}{y} \tilde{c}(t, x-y) \tilde{c}(t,y)$$
$$- \sum_{y=1}^{\infty} \frac{K(x,y)}{y} \tilde{c}(t,x) \tilde{c}(t,y),$$

where $\tilde{c}(t,x) = x c(t,x)$. Equation (1.5) describes the time evolution of the particle mass density (relative mass of particles of a given size). A corresponding mass flow process is represented by a system of particles

(1.6) $\qquad \tilde{x}_i^{(n)}(t), \qquad i = 1, \ldots, \tilde{N}^{(n)}(t), \qquad \tilde{N}^{(n)}(0) = n,$

with an appropriately modified evolution rule so that the solution to equation (1.5) is approximated. Due to the equivalence of the equations, the mass flow model provides an alternative approximation to the solution of the original coagulation equation (1.1), namely,

$$\sum_{x=1}^{\infty} \varphi(x) c(t,x) = \sum_{x=1}^{\infty} \frac{\varphi(x)}{x} \tilde{c}(t,x)$$

(1.7)

$$= \lim_{n \to \infty} \frac{1}{n} \sum_{i=1}^{\tilde{N}^{(n)}(t)} \frac{\varphi(\tilde{x}_i^{(n)}(t))}{\tilde{x}_i^{(n)}}, \qquad t \geq 0.$$

We refer to [9] concerning a rigorous proof. A time-discrete version of the mass flow process was introduced and studied in [2].

In the direct simulation model (1.3), the mass of the system is conserved so that (1.4) implies

$$\sum_{x=1}^{\infty} x c(t,x) \sim \frac{1}{n} \sum_{i=1}^{N^{(n)}(t)} x_i^{(n)}(t) = \frac{1}{n} \sum_{i=1}^{n} x_i^{(n)}(0), \qquad t \geq 0.$$

Thus, there is no way to approximate the gelation point following the definition (1.2). Instead, the gelation effect is related to the formation of a very big cluster (order of the whole system) in finite time. In the mass flow model (1.6), the situation is different. According to (1.7), one obtains

$$\sum_{x=1}^{\infty} x c(t,x) \sim \frac{\tilde{N}^{(n)}(t)}{n}, \qquad t \geq 0,$$

that is, the mass of the solution to equation (1.1) is approximated by the normalized number of particles in the mass flow system. Since the particles grow, there is a chance that their number drops due to some explosion phenomenon (infinitely many jumps in finite time). It was conjectured in [9] that, in case of gelling kernels, the mass flow process explodes and that the gelation time (1.2) is the limit (as $n \to \infty$) of the (random) explosion times of the approximating finite particle systems.

The paper contains three main results. *First* some rather general explosion criteria for pure jump processes are obtained. Using ideas from [18, 19], previously known results for one-dimensional processes are generalized to processes with a locally compact state space. This allows one to treat many stochastic coagulation–fragmentation models. *Second* we prove that the mass flow process explodes almost surely for a very wide class of gelling kernels in the case of pure coagulation. This result confirms the first part of the conjecture mentioned above. Its simple proof illustrates the usefulness of the general explosion criteria. *Third* we study the case of pure multiple fragmentation for a continuous size space. Applying the general criteria, we prove explosion results for both the direct simulation model and the mass flow model, when the fragmentation rate at zero grows sufficiently fast. These results correspond to another kind of phase transition—a loss of mass to zero (while gelation corresponds to a loss of mass to infinity). This phenomenon was studied a long time ago in [12] and named "transformation into dust." Later it was called "shattering" in the physical literature [23, 30].

The paper is organized as follows. In Section 2 we introduce the minimal jump process with an arbitrary locally compact separable metric state space and prove several sufficient criteria for explosion of this process. We give some examples and reproduce the previously known conditions for the one-dimensional case. In Section 3 we define the direct simulation model and the mass flow model. These models are extended to include both coagulation and multiple fragmentation, as well as source and efflux terms, which are important in many applications. The relation of the stochastic models to the corresponding macroscopic equations is sketched heuristically, referring to convergence results available in the literature. In Section 4 we apply the explosion criteria from Section 2 to the stochastic coagulation–fragmentation models introduced in Section 3. First we consider the direct simulation model and prove a simple nonexplosion criterion in the general case, as well as an explosion result in the pure fragmentation case. Then we consider the mass flow model. In the pure coagulation case, we prove that there is almost sure explosion for arbitrary homogeneous coagulation kernels with exponent bigger than 1. Finally, we derive an explosion result in the pure fragmentation case and show that the explosion properties of both models are not equivalent. Section 5 contains some concluding remarks.

## 2. Explosion criteria for jump processes.

2.1. *The minimal jump process.* Let $E$ and $E'$ be separable metric spaces. The sets of measurable and continuous functions on $E$ are denoted by $M(E)$ and $C(E)$, respectively. Furthermore, $\mathcal{M}_b(E)$ and $\mathcal{P}(E)$ are the sets of bounded Borel measures and of probability measures on the Borel-$\sigma$-algebra $\mathcal{B}(E)$. Finally, let $\mathbb{1}_B$ denote the indicator function of a set $B$, and $\delta_\xi$ be the Dirac measure on $\xi \in E$. A kernel from $E$ to $E'$ (on $E$ if $E = E'$) is a function $q : E \times \mathcal{B}(E') \to [0, \infty)$ such that

$$q(\cdot, B) \in M(E) \quad \forall B \in \mathcal{B}(E') \quad \text{and} \quad q(\xi, \cdot) \in \mathcal{M}_b(E') \quad \forall \xi \in E.$$

A kernel $q$ is called compactly bounded if

$$\sup_{\xi \in C} q(\xi, E') < \infty \quad \text{for any compact } C \subset E.$$

Let $q$ be a compactly bounded kernel on a locally compact separable metric space $E$. Let $\zeta_0, \zeta_1, \ldots$ be a Markov chain in $E$ with initial distribution $\nu_0 \in \mathcal{P}(E)$ and transition function $p : E \times \mathcal{B}(E) \to [0, 1]$ defined by

$$(2.1) \qquad p(\xi, B) = \begin{cases} \dfrac{q(\xi, B)}{\lambda(\xi)} : \lambda(\xi) > 0, \\ \mathbb{1}_B(\xi) : \lambda(\xi) = 0, \end{cases}$$

where

(2.2) $$\lambda(\xi) = q(\xi, E), \qquad \xi \in E.$$

Let $T_0, T_1, \ldots$ be independent and exponentially distributed random variables with mean 1 that are also independent of $(\zeta_k)$, all defined on some probability space $(\Omega, \mathcal{F}, \mathbb{P})$. Introduce the jump times

$$\tau_0 = 0, \qquad \tau_l = \sum_{k=0}^{l-1} \frac{T_k}{\lambda(\zeta_k)}, \qquad l = 1, 2, \ldots,$$

where $T_k/0 := \infty$, and the explosion time

(2.3) $$\tau_\infty = \lim_{l \to \infty} \tau_l = \sum_{k=0}^{\infty} \frac{T_k}{\lambda(\zeta_k)}.$$

The *minimal jump process*, corresponding to the kernel $q$ and the initial distribution $\nu_0$, is defined as (cf. [11], page 263 and [25], page 69)

(2.4) $$\zeta^\Delta(t) = \begin{cases} \zeta_l : \tau_l \leq t < \tau_{l+1}, \\ \Delta : t \geq \tau_\infty, \end{cases} \qquad t \geq 0,$$

where $\Delta \notin E$ determines the one-point compactification of $E$ (cf. [3], page 205). The process is called *regular*, if

(2.5) $$\mathbb{P}(\tau_\infty = \infty) = 1.$$

Otherwise, the process is called *explosive*. Concerning the history of the subject we refer to [6], Prologue.

Regularity and explosion of the minimal jump process can be studied using the properties (cf. [25], page 71)

(2.6) $$\sum_{k=0}^{\infty} \frac{1}{a_k} < \infty \quad \Longrightarrow \quad \mathbb{P}\left(\sum_{k=0}^{\infty} \frac{T_k}{a_k} < \infty\right) = 1$$

and

(2.7) $$\sum_{k=0}^{\infty} \frac{1}{a_k} = \infty \quad \Longrightarrow \quad \mathbb{P}\left(\sum_{k=0}^{\infty} \frac{T_k}{a_k} = \infty\right) = 1,$$

for any nonnegative sequence $(a_k)$. The independence of $(T_k)$ and $(\zeta_k)$ allows one to conclude from (2.6) and (2.7) that [cf. (2.3)]

$$\mathbb{P}\left(\tau_\infty < \infty \Big| \sum_{k=0}^{\infty} \frac{1}{\lambda(\zeta_k)} = \infty\right) = \mathbb{P}\left(\tau_\infty = \infty \Big| \sum_{k=0}^{\infty} \frac{1}{\lambda(\zeta_k)} < \infty\right) = 0,$$

which implies

$$\{\tau_\infty < \infty\} = \left\{\sum_{k=0}^{\infty} \frac{1}{\lambda(\zeta_k)} < \infty\right\} \qquad \text{almost surely.}$$

In particular, a necessary and sufficient condition for regularity (2.5) is (cf. [5], page 337)

$$\mathbb{P}\left(\sum_{k=0}^{\infty}\frac{1}{\lambda(\zeta_k)}=\infty\right)=1.$$

Correspondingly, the process is explosive if and only if

(2.8) $$\mathbb{P}\left(\sum_{k=0}^{\infty}\frac{1}{\lambda(\zeta_k)}<\infty\right)>0.$$

Note that boundedness of the waiting time parameter $\lambda$ implies regularity, and that

(2.9) $$\sum_{k=0}^{\infty}\frac{1}{\lambda(\zeta_k)}<\infty \quad \Longrightarrow \quad \lim_{k\to\infty}\zeta_k=\Delta.$$

Thus, a necessary condition for explosion is

(2.10) $$\mathbb{P}\left(\lim_{k\to\infty}\zeta_k=\Delta\right)>0.$$

2.2. *Explosion criteria.* First we formulate a general result concerning explosion of the minimal jump process defined in the previous section.

THEOREM 2.1. *Let $q$ be a compactly bounded kernel on a locally compact separable metric space $E$ and $(\zeta_k)$ be the corresponding Markov chain [cf. (2.1) and (2.2)]. Consider the sets*

(2.11) $$E_\varepsilon(\eta) = \left\{\xi \in E : \int_E [\eta(\xi_1)-\eta(\xi)]q(\xi,d\xi_1) \geq \varepsilon\right\}$$

*and*

(2.12) $\Omega_\varepsilon(\eta) = \{\omega \in \Omega : \zeta_k(\omega) \in E_\varepsilon(\eta),\ \forall k \geq \bar{k}(\omega),\ \text{for some } \bar{k}(\omega)\},$

*where $\eta$ is a bounded measurable function on $E$ and $\varepsilon \geq 0$. Then*

(2.13) $$\sum_{k=0}^{\infty}\frac{1}{\lambda(\zeta_k)}<\infty \qquad a.s.\ on\ \Omega_\varepsilon(\eta),\ \forall \varepsilon>0,$$

*that is, the minimal jump process, corresponding to the kernel $q$, explodes almost surely on the set of all trajectories staying in $E_\varepsilon(\eta)$ for sufficiently large $k$ and some $\varepsilon > 0$.*

REMARK 2.2. According to Theorem 2.1, a sufficient condition for explosion is the existence of some bounded measurable function $\eta$ such that

(2.14) $$\mathbb{P}(\Omega_\varepsilon(\eta))>0 \qquad \text{for some } \varepsilon>0.$$

It follows from (2.9) and (2.13) that

$$\lim_{k\to\infty} \zeta_k = \Delta \qquad \text{a.s. on } \Omega_\varepsilon(\eta),\ \forall\,\varepsilon > 0.$$

Thus, condition (2.14) puts some restriction on the possible structure of the set (2.11). In particular, this set can not be contained in any compact.

Before proving Theorem 2.1, we derive several corollaries providing sufficient conditions for explosion.

COROLLARY 2.3. *Suppose that*

$$(2.15) \qquad \int_E [\eta(\xi_1) - \eta(\xi)] q(\xi, d\xi_1) \geq \varepsilon \qquad \forall\,\xi \in E^+,$$

*for some bounded measurable function $\eta$, some subset $E^+ \subset E$ and some $\varepsilon > 0$. Then*

$$\sum_{k=0}^{\infty} \frac{1}{\lambda(\zeta_k)} < \infty \qquad \text{a.s. on } \{\zeta_k \in E^+,\ \forall\,k\},$$

*that is, the minimal jump process explodes almost surely on the set of all trajectories living in $E^+$.*

PROOF. Assumption (2.15) implies $E^+ \subset E_\varepsilon(\eta)$ and $\{\zeta_k \in E^+,\ \forall\,k\} \subset \Omega_\varepsilon(\eta)$ so that the assertion follows from (2.13). □

A particularly simple criterion for explosion is obtained in the case $E^+ = E$.

COROLLARY 2.4. *Suppose there exists some bounded measurable function $\eta$ such that*

$$(2.16) \qquad \int_E [\eta(\xi_1) - \eta(\xi)] q(\xi, d\xi_1) \geq \varepsilon \qquad \forall\,\xi \in E,$$

*for some $\varepsilon > 0$. Then*

$$\mathbb{P}\bigg(\sum_{k=0}^{\infty} \frac{1}{\lambda(\zeta_k)} < \infty\bigg) = 1,$$

*that is, the minimal jump process explodes almost surely, for any initial distribution.*

REMARK 2.5. Under the assumptions of Corollary 2.4, it follows from (2.9) that

$$\mathbb{P}\Big(\lim_{k\to\infty} \zeta_k = \Delta\Big) = 1.$$

COROLLARY 2.6. *Suppose there exists a bounded measurable function $\eta$ such that*

$$\int_E [\eta(\xi_1) - \eta(\xi)] q(\xi, d\xi_1) \geq \varepsilon \qquad \forall \xi \notin C, \tag{2.17}$$

*for some compact $C \subset E$ and some $\varepsilon > 0$. Then*

$$\sum_{k=0}^{\infty} \frac{1}{\lambda(\zeta_k)} < \infty \qquad \text{a.s. on } \left\{ \lim_{k \to \infty} \zeta_k = \Delta \right\}, \tag{2.18}$$

*that is, the minimal jump process explodes almost surely on the set of all trajectories tending to the compactification point.*

PROOF. Assumption (2.17) implies $E \setminus C \subset E_\varepsilon(\eta)$ and $\{\lim_{k \to \infty} \zeta_k = \Delta\} \subset \Omega_\varepsilon(\eta)$ so that the assertion follows from (2.13). □

REMARK 2.7. Under the assumptions of Corollary 2.6, it follows from (2.9) and (2.18) that

$$\left\{ \sum_{k=0}^{\infty} \frac{1}{\lambda(\zeta_k)} < \infty \right\} = \left\{ \lim_{k \to \infty} \zeta_k = \Delta \right\} \qquad \text{a.s.}$$

so that the necessary condition (2.10) is also sufficient for explosion of the minimal jump process. Note that trajectories may stay in the compact set $C$ forever.

The proof of Theorem 2.1 is prepared by the following lemma.

LEMMA 2.8. *Let $\eta$ be a bounded measurable function on $E$. Then*

$$\exists \sum_{k=0}^{\infty} [\mathbb{E}(\eta(\zeta_{k+1})|\zeta_k) - \eta(\zeta_k)] \qquad \text{finite a.s. on } \Omega_0(\eta),$$

*that is, the infinite sum is finite for almost all $\omega \in \Omega_0(\eta)$ [cf. (2.12)].*

PROOF. The sequence

$$W_n = \sum_{k=0}^{n-1} [\mathbb{E}(\eta(\zeta_{k+1})|\zeta_k) - \eta(\zeta_k)] - \eta(\zeta_n), \qquad n \geq 1, W_0 = -\eta(\zeta_0), \tag{2.19}$$

is a martingale with respect to the filtration of $(\zeta_k)$. Representing it in the form

$$W_n = \sum_{k=0,\ldots,n-1 \,:\, \zeta_k \in E_0(\eta)} [\mathbb{E}(\eta(\zeta_{k+1})|\zeta_k) - \eta(\zeta_k)]$$

$$+ \sum_{k=0,\ldots,n-1 \,:\, \zeta_k \notin E_0(\eta)} [\mathbb{E}(\eta(\zeta_{k+1})|\zeta_k) - \eta(\zeta_k)] - \eta(\zeta_n) \tag{2.20}$$

and introducing the sequence of stopping times

$$\sigma_N = \inf\{k > \sigma_{N-1} : \zeta_k \notin E_0(\eta)\}, \qquad N \geq 2,$$
$$\sigma_1 = \inf\{k \geq 0 : \zeta_k \notin E_0(\eta)\},$$

one concludes that

$$W_{\min(n,\sigma_N)} \geq -\sup_{\xi \in E}|\eta(\xi)|[2\#\{k=0,\ldots,\min(n,\sigma_N)-1 : \zeta_k \notin E_0(\eta)\}+1]$$
$$\geq -(2N-1)\sup_{\xi \in E}|\eta(\xi)| \qquad \forall n, N \geq 1.$$

Since the stopped process $W_{\min(n,\sigma_N)}$ is a martingale and bounded from below, it has a.s. finite limits (cf., e.g., [28], Sections II.49 and 57), that is,

$$\mathbb{P}\bigg(\exists \lim_{n \to \infty} W_{\min(n,\sigma_N)} \text{ finite}\bigg) = 1 \qquad \forall N.$$

Consequently, one obtains

(2.21) $\qquad \exists \lim_{n \to \infty} W_n \qquad$ finite a.s. on $\{\sigma_N = \infty\} \, \forall N.$

Note that the first sum at the right-hand side of (2.20) increases, while the second sum has at most $N-1$ elements (if $\sigma_N = \infty$) and the last term is bounded. Thus, (2.21) implies

(2.22) $\qquad \exists \lim_{n \to \infty} \eta(\zeta_n) \qquad$ finite a.s. on $\{\sigma_N = \infty\} \, \forall N.$

Finally, one concludes from (2.19), (2.21) and (2.22) that

$$\exists \sum_{k=0}^{\infty}[\mathbb{E}(\eta(\zeta_{k+1})|\zeta_k) - \eta(\zeta_k)] \qquad \text{finite a.s. on } \{\sigma_N = \infty\} \, \forall N,$$

and the assertion follows from the fact that $\Omega_0(\eta) = \bigcup_N \{\sigma_N = \infty\}$. $\square$

PROOF OF THEOREM 2.1. Note that

$$\int_E [\eta(\xi_1) - \eta(\xi)]q(\xi, d\xi_1) = \lambda(\xi)[\mathbb{E}(\eta(\zeta_{k+1})|\zeta_k = \xi) - \eta(\xi)]$$

and

$$\frac{1}{\lambda(\zeta_k)} \leq \frac{1}{\varepsilon}[\mathbb{E}(\eta(\zeta_{k+1})|\zeta_k) - \eta(\zeta_k)] \qquad \text{if } \zeta_k \in E_\varepsilon(\eta).$$

Since

$$\sum_{k=0}^{\infty} \frac{1}{\lambda(\zeta_k)} \leq \sum_{k=0}^{\bar{k}-1} \frac{1}{\lambda(\zeta_k)} + \frac{1}{\varepsilon}\sum_{k=\bar{k}}^{\infty}[\mathbb{E}(\eta(\zeta_{k+1})|\zeta_k) - \eta(\zeta_k)] \qquad \text{on } \Omega_\varepsilon(\eta),$$

the assertion follows from Lemma 2.8. $\square$

2.3. *Examples.*

EXAMPLE 2.9. Consider the case
$$E = \{1, 2, \dots\}, \qquad q(\xi, d\xi_1) = \lambda(\xi)\delta_{\xi+1}(d\xi_1), \qquad \zeta_0 = 1.$$
The trajectory of the Markov chain is deterministic, $\zeta_k = k+1, k \geq 0$, and the necessary and sufficient condition for explosion (2.8) takes the form

(2.23) $$\sum_{k=1}^{\infty} \frac{1}{\lambda(k)} < \infty.$$

First we show that the sufficient condition (2.16) of Corollary 2.4 can always be satisfied. Indeed, choosing the function
$$\eta(\xi) = \sum_{k=1}^{\xi-1} \frac{1}{\lambda(k)}, \qquad \xi \geq 2, \ \eta(1) = 0,$$
which is bounded if (2.23) holds, one obtains
$$[\eta(\xi+1) - \eta(\xi)]\lambda(\xi) = 1 \qquad \forall \xi \in E.$$
Furthermore, for $\lambda(\xi) = \xi$ and any strictly increasing bounded function $\eta$, this example illustrates that the condition
$$\int_E [\eta(\xi_1) - \eta(\xi)]q(\xi, d\xi_1) > 0 \qquad \forall \xi \in E,$$
instead of (2.16), would not be sufficient. Finally, for $\eta(\xi) = \xi$ and constant $\lambda(\xi)$, this example illustrates that condition (2.16) for unbounded $\eta$ would not be sufficient.

EXAMPLE 2.10. Consider the one-dimensional case $E = [1, \infty)$ and the bounded measurable function
$$\eta(\xi) = -\xi^{-\alpha}, \qquad \xi \in E \text{ for some } \alpha > 0.$$
Condition (2.17) takes the form

(2.24) $$\lambda(\xi)\xi^{-\alpha} \mathbb{E}\left(1 - \left(\frac{\xi}{\zeta_1}\right)^{\alpha} \Big| \zeta_0 = \xi\right) \geq \varepsilon \qquad \text{for sufficiently large } \xi.$$

With the notation
$$n_{-\alpha}(\xi) = \mathbb{E}\left(1 - \left(\frac{\xi}{\zeta_1}\right)^{\alpha} \Big| \zeta_0 = \xi\right) \quad \text{and} \quad f(\xi) = \varepsilon \xi^{\alpha+1},$$
condition (2.24) transforms into
$$\lambda(\xi)\xi n_{-\alpha}(\xi) \geq f(\xi) \qquad \text{for sufficiently large } \xi,$$
so that the result of [18], Theorem 3 is basically reproduced by Corollary 2.6.

**3. Deterministic and stochastic coagulation–fragmentation models.** In this section we introduce two stochastic coagulation–fragmentation models that are covered by the framework of Section 2. Explosion phenomena in these models will be studied in Section 4.

The two models are related to each other in the sense of (1.4) and (1.7). In order to demonstrate this, we let them depend on a parameter $n = 1, 2, \ldots$ and illustrate heuristically their relation (when $n \to \infty$) to corresponding macroscopic equations of the type (1.1) and (1.5). Rigorous convergence results and more references can be found, for example, in [10]. The treatment of the macroscopic equations in the general setting (including coagulation, multiple fragmentation, source and efflux terms) is of independent interest, since many equations known from the literature (weak and strong, discrete and continuous) are covered in a unified way.

We consider stochastic processes with the state space

$$(3.1) \qquad E^{(n)} = \left\{ \frac{1}{n} \sum_{i=1}^{N} \delta_{x_i} : N \geq 0, \ x_i \in X, \ i = 1, \ldots, N \right\},$$

where $X$ is a locally compact separable metric space. A state $\xi \in E^{(n)}$ is interpreted as a system of particles with types from $X$ and weights $\frac{1}{n}$. Standard examples are

$$(3.2) \qquad X = \{1, 2, \ldots\} \quad \text{or} \quad X = (0, \infty),$$

corresponding to discrete and continuous particles sizes, respectively. However, many of the constructions of this section are valid in the general setup, which allows one to cover particles consisting of several species, and other cases of practical interest.

REMARK 3.1. The space $E^{(n)}$, endowed with an appropriate metric of weak convergence, is locally compact and separable (cf. Remark 2.10 and Lemma 5.1 in [10]). Note that $\xi_k \to \Delta$ in $E^{(n)}$ if either $N_k \to \infty$ or $x_{k,i} \to \Delta_X$ for some $i$, where $\Delta_X$ denotes the compactification point of $X$.

3.1. *The direct simulation model.* Consider the kernel

$$(3.3) \quad \begin{aligned} q^{(n)}&(\xi, B) \\ &= n \int_X \mathbb{1}_B(J_S(\xi, x)) S(dx) \\ &\quad + \sum_{i=1}^N \mathbb{1}_B(J_e(\xi, i)) e(x_i) + \sum_{i=1}^N \int_Z \mathbb{1}_B(J_F(\xi, i, z)) F(x_i, dz) \\ &\quad + \frac{1}{2n} \sum_{1 \leq i \neq j \leq N} \mathbb{1}_B(J_K(\xi, i, j)) K(x_i, x_j), \qquad \xi \in E^{(n)}, \ B \in \mathcal{B}(E^{(n)}), \end{aligned}$$

with the jump transformations

$$J_S(\xi, x) = \xi + \frac{1}{n}\delta_x,$$

$$J_e(\xi, i) = \xi - \frac{1}{n}\delta_{x_i},$$

(3.4)

$$J_F(\xi, i, z) = \xi - \frac{1}{n}\delta_{x_i} + \frac{1}{n}[\delta_{z_1} + \cdots + \delta_{z_k}],$$

$$J_K(\xi, i, j) = \xi - \frac{1}{n}[\delta_{x_i} + \delta_{x_j}] + \frac{1}{n}\delta_{x_i+x_j}.$$

Here $S \in \mathcal{M}_b(X)$ is some source measure. The nonnegative function $e \in C(X)$ denotes the efflux intensity. The compactly bounded fragmentation kernel $F$ from $X$ to

(3.5) $$Z = \bigcup_{k=2}^{\infty} X^k$$

is assumed to satisfy the mass conservation property

(3.6) $$F(x, Z \setminus Z(x)) = 0 \qquad \forall\, x \in X,$$

where $Z(x) = \{z \in Z : z_1 + \cdots + z_k = x\}$. The nonnegative function $K \in C(X \times X)$ denotes the coagulation kernel. Under the above assumptions, the kernel (3.3) is compactly bounded (cf. [10], Lemma 5.1). The corresponding minimal jump process is called *direct simulation process*.

Note that

(3.7)
$$\lambda^{(n)}(\xi) = q^{(n)}(\xi, E^{(n)})$$
$$= n\bigg[S(X) + \int_X e(x)\xi(dx)$$
$$+ \int_X F(x, Z)\xi(dx) + \frac{1}{2n^2}\sum_{1 \le i \ne j \le N} K(x_i, x_j)\bigg]$$

and

(3.8)
$$\int_{E^{(n)}}\bigg[\int_X \varphi(x)\xi_1(dx) - \int_X \varphi(x)\xi(dx)\bigg] q^{(n)}(\xi, d\xi_1)$$
$$= \int_X \varphi(x) S(dx) - \int_X \varphi(x) e(x) \xi(dx)$$
$$+ \int_X \int_Z [\varphi(z_1) + \cdots + \varphi(z_k) - \varphi(x)] F(x, dz)\xi(dx)$$
$$+ \frac{1}{2n^2} \sum_{1 \le i \ne j \le N} [\varphi(x_i + x_j) - \varphi(x_i) - \varphi(x_j)] K(x_i, x_j),$$

for any $\xi \in E^{(n)}$ and appropriate test functions $\varphi$. Dynkin's formula and (3.8) suggest that the direct simulation process corresponds (for $n \to \infty$) to the macroscopic *coagulation–fragmentation equation*

$$
\begin{aligned}
\frac{d}{dt} &\int_X \varphi(x)\mu(t,dx) \\
&= \int_X \varphi(x)S(dx) - \int_X \varphi(x)e(x)\mu(t,dx) \\
&\quad + \int_X \int_Z [\varphi(z_1) + \cdots + \varphi(z_k) - \varphi(x)]F(x,dz)\mu(t,dx) \\
&\quad + \frac{1}{2}\int_X \int_X [\varphi(x+y) - \varphi(x) - \varphi(y)] \\
&\quad\quad \times K(x,y)\mu(t,dx)\mu(t,dy), \qquad t > 0,
\end{aligned}
\tag{3.9}
$$

with the initial condition $\mu(0,dx) = \mu_0(dx)$, where $\mu_0$ is the limit of the initial state of the process.

Let $F^{(k)}$ denote the restrictions of the fragmentation kernel $F$ to the sets $X^k$ [cf. (3.5)]. Define the *symmetrized fragmentation kernel* $F_{\text{sym}}$ by its corresponding restrictions

$$
\begin{aligned}
F^{(k)}_{\text{sym}}&(x,dz_1,\ldots,dz_k) \\
&= \frac{1}{k!}\sum_\pi F^{(k)}(x,dz_{\pi(1)},\ldots,dz_{\pi(k)}), \qquad k = 2,3,\ldots,
\end{aligned}
\tag{3.10}
$$

where the sum is taken over all permutations of $\{1,\ldots,k\}$. Introduce the 1-marginals of the kernels (3.10)

$$
F^{(k|1)}_{\text{sym}}(x,dy) = F^{(k)}_{\text{sym}}(x,dy,X,\ldots,X)
\tag{3.11}
$$

and the kernel

$$
F^{(1)}_{\text{sym}}(x,dy) = \sum_{k=2}^\infty k F^{(k|1)}_{\text{sym}}(x,dy).
\tag{3.12}
$$

Then one obtains

$$
\begin{aligned}
\int_Z &[\varphi(z_1) + \cdots + \varphi(z_k)]F(x,dz) \\
&= \sum_{k=2}^\infty \int_{X^k} [\varphi(z_1) + \cdots + \varphi(z_k)]F^{(k)}_{\text{sym}}(x,dz_1,\ldots,dz_k) \\
&= \sum_{k=2}^\infty k \int_X \varphi(y) F^{(k|1)}_{\text{sym}}(x,dy) = \int_X \varphi(y) F^{(1)}_{\text{sym}}(x,dy).
\end{aligned}
\tag{3.13}
$$

Using the *symmetrized coagulation kernel*

(3.14) $$K_{\text{sym}}(x,y) = \tfrac{1}{2}[K(x,y) + K(y,x)],$$

one obtains

(3.15)
$$\begin{aligned}&\tfrac{1}{2}\int_X\int_X[\varphi(x+y)-\varphi(x)-\varphi(y)]K(x,y)\mu(t,dx)\mu(t,dy)\\&=\tfrac{1}{2}\int_X\int_X[\varphi(x+y)-\varphi(x)-\varphi(y)]K_{\text{sym}}(x,y)\mu(t,dx)\mu(t,dy)\\&=\int_X\int_X[\tfrac{1}{2}\varphi(x+y)-\varphi(x)]K_{\text{sym}}(x,y)\mu(t,dx)\mu(t,dy).\end{aligned}$$

REMARK 3.2. According to (3.13) and (3.15), the macroscopic equation (3.9) does not change if the fragmentation and coagulation kernels of the direct simulation process are replaced by their symmetrizations.

Preparing the transition from (3.9) to an equation for the densities

$$\mu(t,dx) = c(t,x)\,dx,$$

we assume in the continuous case $X = (0,\infty)$ that

(3.16) $$S(dx) = s(x)\,dx$$

and [cf. (3.12)]

(3.17) $$F_{\text{sym}}^{(1)}(x,dy) = f_{\text{sym}}^{(1)}(x,y)\,dy.$$

Note the identity

(3.18) $$\int_0^\infty\int_0^\infty \psi(x,y)\,dy\,dx = \int_0^\infty\int_0^x \psi(x-y,y)\,dy\,dx,$$

where $\psi$ is an appropriate test function. Using (3.13), (3.15) and (3.18), one obtains the *continuous coagulation–fragmentation equation* with source and efflux terms

(3.19)
$$\begin{aligned}\frac{\partial}{\partial t}c(t,x) &= s(x) - e(x)c(t,x)\\&\quad + \int_0^\infty f_{\text{sym}}^{(1)}(x+y,x)c(t,x+y)\,dy - c(t,x)F(x,Z)\\&\quad + \frac{1}{2}\int_0^x K_{\text{sym}}(x-y,y)c(t,x-y)c(t,y)\,dy\\&\quad - \int_0^\infty K_{\text{sym}}(x,y)c(t,x)c(t,y)\,dy.\end{aligned}$$

In the discrete case $X = \{1,2,\dots\}$, an equation analogous to (3.19) is obtained, with integrals replaced by sums.

REMARK 3.3. According to (3.13), the kernel (3.12) can be expressed in the form

$$F^{(1)}_{\text{sym}}(x, dy) = \int_Z [\delta_{z_1}(dy) + \cdots + \delta_{z_k}(dy)] F(x, dz). \tag{3.20}$$

Thus, the quantity

$$\frac{1}{F(x,Z)} F^{(1)}_{\text{sym}}(x, (0,y)) = \mathbb{E}_x \sum_{i=1}^k \mathbb{1}_{(0,y)}(z_i) = \mathbb{E}_x \sum_{i\,:\,z_i<y} 1$$

represents the average number of fragments with size less than $y$, resulting from the fragmentation of a particle of size $x$. In particular, the average number of fragments is expressed as

$$\frac{1}{F(x,Z)} F^{(1)}_{\text{sym}}(x, X). \tag{3.21}$$

Note that $F^{(1)}_{\text{sym}}(x, [\varepsilon, x)) < \infty, \forall \varepsilon \in (0, x)$, due to mass conservation (3.6). Thus, for test functions $\varphi$ with compact support, the corresponding term in (3.9) is finite even if the average number of fragments is infinite. Finally, we mention that mass conservation implies

$$x F^{(k)}(x, X^k) = \int_{X^k} (z_1 + \cdots + z_k) F^{(k)}(x, dz_1, \ldots, dz_k) = k \int_X y F^{(k|1)}_{\text{sym}}(x, dy)$$

and

$$F(x, Z) = \frac{1}{x} \int_X y F^{(1)}_{\text{sym}}(x, dy). \tag{3.22}$$

EXAMPLE 3.4. Consider the *binary fragmentation case*

$$F(x, dz) = F^{(2)}(x, dz_1, dz_2) = \tfrac{1}{2} F^{(1)}(x, dz_1) \delta_{x-z_1}(dz_2), \tag{3.23}$$

where $F^{(1)}$ is a kernel on $X$ such that

$$F^{(1)}(x, X \setminus (0, x)) = 0 \qquad \forall\, x \in X.$$

One obtains [cf. (3.10)–(3.12)]

$$F^{(2)}_{\text{sym}}(x, dz_1, dz_2) = \tfrac{1}{4}[F^{(1)}(x, dz_1) \delta_{x-z_1}(dz_2) + F^{(1)}(x, dz_2) \delta_{x-z_2}(dz_1)]$$

so that

$$F^{(2|1)}_{\text{sym}}(x, dy) = \tfrac{1}{4}\left[F^{(1)}(x, dy) + \int_X F^{(1)}(x, dx_1) \delta_{x-x_1}(dy)\right]$$

and

$$F^{(1)}_{\text{sym}}(x, dy) = \tfrac{1}{2}\left[F^{(1)}(x, dy) + \int_X F^{(1)}(x, dx_1) \delta_{x-x_1}(dy)\right]. \tag{3.24}$$

Assuming
$$F^{(1)}(x, dy) = \mathbb{1}_{(0,x)}(y) f^{(1)}(x, y)\, dy,$$

one obtains from (3.17) and (3.24) that
$$f^{(1)}_{\text{sym}}(x, y) = \tfrac{1}{2}[f^{(1)}(x, y) + f^{(1)}(x, x - y)].$$

Since
$$F(x, Z) = \tfrac{1}{2}\int_0^x f^{(1)}(x, y)\, dy = \tfrac{1}{2}\int_0^x f^{(1)}_{\text{sym}}(x, y)\, dy,$$

the fragmentation term in (3.19) takes the usual form. Note that (3.22) implies
$$\frac{1}{2} F^{(1)}(x, X) = \frac{1}{x}\int_X y F^{(1)}_{\text{sym}}(x, dy) \qquad \forall\, x \in X.$$

3.2. *The mass flow model.* The solution $\mu(t, dx)$ of equation (3.9) represents the flow of concentration in the size space $X$. Since the total mass of the system is determined as $\int_X x\mu(t, dx)$, we call the function
$$\tilde{\mu}(t, dx) = x\mu(t, dx), \qquad t \geq 0,$$

the mass flow. Considering test functions of the form $\varphi(x) = x\psi(x)$, one obtains [cf. (3.9) and (3.13)]
$$\int_X \int_X y\psi(y) F^{(1)}_{\text{sym}}(x, dy)\mu(t, dx) = \int_X \int_X \psi(y) \tilde{F}(x, dy)\tilde{\mu}(t, dx),$$

where

(3.25) $$\tilde{F}(x, dy) = \frac{y}{x} F^{(1)}_{\text{sym}}(x, dy),$$

and [cf. (3.15)]
$$\int_X \int_X \left[\frac{1}{2}(x+y)\psi(x+y) - x\psi(x)\right] K_{\text{sym}}(x, y)\mu(t, dx)\mu(t, dy)$$
$$= \int_X \int_X [x\psi(x+y) - x\psi(x)] K_{\text{sym}}(x, y)\mu(t, dx)\mu(t, dy)$$
$$= \int_X \int_X [\psi(x+y) - \psi(x)] \frac{K_{\text{sym}}(x, y)}{y} \tilde{\mu}(t, dx)\tilde{\mu}(t, dy).$$

Since, according to (3.22),

(3.26) $$\tilde{F}(x, X) = F(x, Z),$$

equation (3.9) takes the form of the *mass flow equation*

$$\begin{aligned}(3.27)\quad &\frac{d}{dt}\int_X \psi(x)\tilde{\mu}(t,dx)\\ &= \int_X \psi(x)xS(dx)\\ &\quad - \int_X \psi(x)e(x)\tilde{\mu}(t,dx) + \int_X\int_X [\psi(y)-\psi(x)]\tilde{F}(x,dy)\tilde{\mu}(t,dx)\\ &\quad + \int_X\int_X [\psi(x+y)-\psi(x)]\frac{K_{\text{sym}}(x,y)}{y}\tilde{\mu}(t,dy)\tilde{\mu}(t,dx),\qquad t>0,\end{aligned}$$

for appropriate test functions $\psi$ and some initial condition.

In the continuous case $X=(0,\infty)$, we assume (3.16) and (3.17). It follows from (3.25) that $\tilde{F}(x,dy)=\tilde{f}(x,y)\,dy$, where

$$\tilde{f}(x,y) = \frac{y}{x}f_{\text{sym}}^{(1)}(x,y).$$

Using (3.18), one obtains from (3.27) an equation for the densities

$$\tilde{\mu}(t,dx) = \tilde{c}(t,x)\,dx,$$

namely, the *continuous mass flow equation* with source and efflux terms

$$\begin{aligned}\frac{\partial}{\partial t}\tilde{c}(t,x) &= xs(x) - e(x)\tilde{c}(t,x)\\ &\quad + \int_0^x \frac{K_{\text{sym}}(x-y,y)}{y}\tilde{c}(t,x-y)\tilde{c}(t,y)\,dy\\ &\quad - \int_0^\infty \frac{K_{\text{sym}}(x,y)}{y}\tilde{c}(t,x)\tilde{c}(t,y)\,dy\\ &\quad + \int_0^\infty \tilde{f}(x+y,x)\tilde{c}(t,x+y)\,dy - \int_0^x \tilde{f}(x,y)\tilde{c}(t,x)\,dy.\end{aligned}$$

This equation is equivalent to (3.19).

We assume

$$(3.28)\qquad \int_X xS(dx) < \infty$$

and introduce the modified kernel [cf. (3.3)]

$$\begin{aligned}(3.29)\quad \tilde{q}^{(n)}(\xi,B) &= n\int_X \mathbb{1}_B(J_S(\xi,x))xS(dx)\\ &\quad + \sum_{i=1}^N \mathbb{1}_B(J_e(\xi,i))e(x_i) + \sum_{i=1}^N \int_X \mathbb{1}_B(\tilde{J}_F(\xi,i,y))\tilde{F}(x_i,dy)\\ &\quad + \frac{1}{n}\sum_{i,j=1}^N \mathbb{1}_B(\tilde{J}_K(\xi,i,j))\frac{K_{\text{sym}}(x_i,x_j)}{x_j},\qquad B\in B(E^{(n)}),\end{aligned}$$

with the modified jump transformations [cf. (3.4)]

$$\tilde{J}_F(\xi, i, y) = \xi - \frac{1}{n}\delta_{x_i} + \frac{1}{n}\delta_y,$$
(3.30)
$$\tilde{J}_K(\xi, i, j) = \xi - \frac{1}{n}\delta_{x_i} + \frac{1}{n}\delta_{x_i+x_j}.$$

The minimal jump process, corresponding to the kernel (3.29), is called *mass flow process*.

Following (3.20), the mass flow fragmentation kernel (3.25) can be represented in the form

$$(3.31) \qquad \tilde{F}(x, dy) = \int_Z \left[\frac{z_1}{x}\delta_{z_1}(dy) + \cdots + \frac{z_k}{x}\delta_{z_k}(dy)\right] F(x, dz).$$

Since the total fragmentation rate does not change [cf. (3.26)], the "fragment" (next state) in the mass flow process is chosen as follows: first fragments $z_1, \ldots, z_k$ are generated according to the normalized direct simulation fragmentation kernel $F(x, dz)$; then one of them is chosen with probabilities proportional to their individual masses.

Note that [cf. (3.7), (3.8) and (3.28)–(3.30)]

$$\tilde{\lambda}^{(n)}(\xi) = \tilde{q}^{(n)}(\xi, E^{(n)})$$
(3.32)
$$= n\bigg[\int_X xS(dx) + \int_X e(x)\xi(dx)$$
$$+ \int_X \tilde{F}(x, X)\xi(dx) + \int_X \int_X \frac{K_{\text{sym}}(x, y)}{y}\xi(dx)\xi(dy)\bigg]$$

and

$$\int_{E^{(n)}} \bigg[\int_X \psi(x)\xi_1(dx) - \int_X \psi(x)\xi(dx)\bigg] \tilde{q}^{(n)}(\xi, d\xi_1)$$
$$= \int_X \psi(x)xS(dx) - \int_X \psi(x)e(x)\xi(dx)$$
(3.33)
$$+ \int_X \int_X [\psi(y) - \psi(x)]\tilde{F}(x, dy)\xi(dx)$$
$$+ \int_X \int_X [\psi(x+y) - \psi(x)]\frac{K_{\text{sym}}(x, y)}{y}\xi(dx)\xi(dy),$$

for any $\xi \in E^{(n)}$ and appropriate test functions $\psi$. The corresponding macroscopic equation (obtained for $n \to \infty$) is (3.27).

3.3. *Comments.* The study of the coagulation–fragmentation equation (3.19) (with $s = e = 0$ and binary fragmentation) goes back to [24].

In the pure fragmentation case $S = 0, e = 0, K = 0$, the mass flow equation (3.27) with test functions $\psi(y) = \mathbb{1}_{(0,x)}(y)$ implies

$$
\begin{aligned}
(3.34) \quad \frac{\partial}{\partial t}\tilde{\mu}(t,(0,x)) &= \int_X \int_X [\mathbb{1}_{(0,x)}(y_1) - \mathbb{1}_{(0,x)}(y)] \tilde{F}(y, dy_1) \tilde{\mu}(t, dy) \\
&= \int_x^\infty \int_0^x \tilde{F}(y, dy_1) \tilde{\mu}(t, dy) \\
&= \int_x^\infty \tilde{F}(y,(0,x)) \, d_y \tilde{\mu}(t,(0,y)).
\end{aligned}
$$

According to (3.31), one obtains

$$
\begin{aligned}
\tilde{F}(y,(0,x)) &= \int_Z \left[ \frac{z_1}{y} \delta_{z_1}((0,x)) + \cdots + \frac{z_k}{y} \delta_{z_k}((0,x)) \right] F(y, dz) \\
&= F(y, Z) \mathbb{E}_y \left[ \frac{1}{y} \sum_{i \,:\, z_i < x} z_i \right]
\end{aligned}
$$

so that equation (3.34) is identical with equation (2) in [12]. Note that $\tilde{\mu}(t,(0,x))$ is the average sum of masses of particles of mass smaller than $x$ at time $t$. Since [cf. (3.26)]

$$
(3.35) \quad \mathbb{E}_y \left[ \frac{1}{y} \sum_{i \,:\, z_i < x} z_i \right] = \frac{\tilde{F}(y,(0,x))}{\tilde{F}(y, X)},
$$

the one-dimensional Markov process introduced in [12], page 279, is just the mass flow process.

**4. Explosion in stochastic coagulation–fragmentation models.** Here we apply the criteria from Section 2 to the models introduced in Section 3. The trajectories of the underlying Markov chain take values in the space (3.1). The state space of single particles is (3.2). The jump kernels and waiting time parameters are determined by (3.3), (3.7) and (3.29), (3.32), respectively. We fix the parameter $n$ and skip the corresponding superscripts.

REMARK 4.1. Condition

$$
(4.1) \quad \mathbb{P}\left( \sum_{k=0}^\infty \frac{1}{\lambda(\zeta_k)} < \infty \Big| \zeta_0 = \xi \right) = 1 \quad \forall \xi \in \hat{E} \subset E,
$$

implies

$$
\mathbb{P}\left( \sum_{k=0}^\infty \frac{1}{\lambda(\zeta_k)} < \infty \right) = \int_{\hat{E}} \mathbb{P}\left( \sum_{k=0}^\infty \frac{1}{\lambda(\zeta_k)} < \infty \Big| \zeta_0 = \xi \right) \nu_0(d\xi) = 1,
$$

that is, explosion with probability one [cf. (2.8)], for any initial distribution $\nu_0$ on $\hat{E}$.

It is often easier to prove (4.1), since trajectories starting at a given $\xi$ remain in a certain part of the state space, where some of the sufficient conditions can be checked. Note that, according to (2.10) and Remark 3.1, explosion implies that (with positive probability) either the number of particles in the system reaches $\infty$ or some of the particles reach 0 or $\infty$.

4.1. *Explosion of the direct simulation process.* We consider the process with the kernel (3.3). In the pure coagulation case the waiting time parameter (3.7) takes the form

$$\lambda(\xi) = \frac{1}{2n} \sum_{1 \leq i \neq j \leq N} K(x_i, x_j).$$

Particle sizes increase so that they cannot reach zero. Both the growth of individual particles and the number of particles in the system are bounded due to mass conservation. So there is no explosion provided the coagulation kernel is bounded on compacts. Analogous arguments apply to the discrete coagulation–fragmentation case (without source term). Thus, all interesting cases in the sense of explosion should include continuous fragmentation. First, we study the situation, when the total fragmentation rate is bounded at zero, and provide some sufficient conditions for regularity. Finally, Theorem 4.3 gives a sufficient condition for explosion in the case of pure fragmentation.

THEOREM 4.2. *Consider the direct simulation kernel* (3.3), *where*

(4.2) $\quad F^{(k)}(x, X^k) = 0 \quad\quad \forall\, x \in X, k \geq k_F,$

(4.3) $\quad\quad F(x, Z) \leq C_F(1+x) \quad\quad \forall\, x \in X,$

(4.4) $\quad\quad K(x,y) \leq C_K(x+y+xy) \quad \forall\, x, y \in X,$

(4.5) $\quad\quad e(x) \leq C_e(1+x) \quad\quad \forall\, x \in X,$

*and*

(4.6) $\quad\quad S(X \cap (C_S, \infty)) = 0,$

*for some constants* $k_F, C_F, C_K, C_e$ *and* $C_S$. *If there is no source term* ($C_S = 0$) *or no coagulation term* ($C_K = 0$), *then the direct simulation process is regular.*

PROOF. According to (4.3)–(4.6), the waiting time parameter (3.7) satisfies

$$\lambda(\xi) \leq n\Big\{ S(X \cap (0, C_S]) + (C_e + C_F)[M_0(\xi) + M_1(\xi)]$$
(4.7)
$$+ \frac{C_K}{2}[2M_0(\xi)M_1(\xi) + M_1(\xi)^2] \Big\},$$

where the notations:

$$M_0(\xi) = \xi(X), \qquad M_1(\xi) = \int_X x\xi(dx), \qquad \xi \in E,$$

are used. Consider trajectories $(\zeta_0, \zeta_1, \dots)$. According to assumption (4.2), the number of particles in the system grows at most linearly, that is,

(4.8) $$\mathbb{P}(M_0(\zeta_k) \leq M_0(\zeta_0) + kk_F, \ \forall k \geq 0) = 1.$$

According to assumption (4.6), the mass of the system grows at most linearly, that is,

(4.9) $$\mathbb{P}(M_1(\zeta_k) \leq M_1(\zeta_0) + kC_S, \ \forall k \geq 0) = 1.$$

If there is no coagulation term, then (4.7)–(4.9) imply (a.s.)

(4.10) $$\lambda(\zeta_k) \leq n\{S(X) + (C_e + C_F)[M_0(\zeta_0) + kk_F + M_1(\zeta_0) + kC_S]\}$$
$$\forall k \geq 0.$$

If there is no source term, then the mass of the system does not grow, that is,

(4.11) $$\mathbb{P}(M_1(\zeta_k) \leq M_1(\zeta_0), \ \forall k \geq 0) = 1,$$

and (4.7), (4.8), (4.11) imply (a.s.)

(4.12) $$\lambda(\zeta_k) \leq n\Big\{(C_e + C_F)[M_0(\zeta_0) + kk_F + M_1(\zeta_0)]$$
$$+ \frac{C_K}{2}[2[M_0(\zeta_0) + kk_F]M_1(\zeta_0) + M_1(\zeta_0)^2]\Big\} \qquad \forall k \geq 0.$$

In both cases (4.10) and (4.12), one obtains (a.s.)

$$\sum_{k=0}^{\infty} \frac{1}{\lambda(\zeta_k)} \geq \sum_{k=0}^{\infty} \frac{1}{C_0 + kC_1} \qquad \text{where } C_0, C_1 < \infty,$$

so that regularity follows. $\square$

THEOREM 4.3. *Consider the direct simulation kernel* (3.3) *in the pure fragmentation case* $(S = 0, e = 0, K = 0)$. *If*

(4.13) $$F(x, Z) \geq \frac{C_F}{x^\alpha} \qquad \forall x \in X, \text{ for some } C_F > 0, \alpha > 0,$$

*then the direct simulation process explodes almost surely, for any initial distribution on the set* $\hat{E} = \{\xi \in E : N \geq 1\}$.

PROOF. In order to make use of Corollary 2.3, we introduce the sets

$$E^+(C) = \left\{\xi \in E : \frac{1}{n}\sum_{i=1}^N x_i \leq C, N \geq 1\right\}, \qquad C > 0, \tag{4.14}$$

and the function $\eta(\xi) = g(\xi(X)), \xi \in E$, where

$$g(x) = \frac{x^\beta}{1+x^\beta}, \qquad x \geq 0, \ 0 < \beta \leq \min(\alpha, 1),$$

is an increasing bounded function. The derivative

$$g'(x) = \frac{\beta x^{\beta-1}}{(1+x^\beta)^2} \tag{4.15}$$

is decreasing so that

$$g(y) - g(x) \geq g'(y)(y-x) \qquad \forall y \geq x \geq 0. \tag{4.16}$$

Furthermore, we note that Jensen's inequality implies (the convex function is $x^{-\alpha}$, and the random variable is uniformly distributed over $N$ different points)

$$\frac{1}{N}\sum_{i=1}^N x_i^{-\alpha} \geq \left(\frac{1}{N}\sum_{i=1}^N x_i\right)^{-\alpha} \qquad \forall x_i > 0, i = 1, \ldots, N, N = 1, 2, \ldots. \tag{4.17}$$

According to (3.3), (3.4), (4.13) and (4.15)–(4.17), one obtains [cf. (4.14)]

$$\begin{aligned}
\int_E [\eta(\xi_1) - \eta(\xi)] q(\xi, d\xi_1) \\
= \sum_{i=1}^N \int_Z [\eta(J_F(\xi, i, z)) - \eta(\xi)] F(x_i, dz) \\
\geq \sum_{i=1}^N \int_Z [g((N+1)/n) - g(N/n)] F(x_i, dz) \\
\geq \frac{C_F}{n} g'((N+1)/n) \sum_{i=1}^N x_i^{-\alpha} \\
\geq \frac{C_F}{n} g'((N+1)/n) N^{1+\alpha} \left(\sum_{i=1}^N x_i\right)^{-\alpha} \\
\geq \frac{C_F}{n} \frac{\beta[(N+1)/n]^{\beta-1}}{(1+[(N+1)/n]^\beta)^2} N^{1+\alpha} (Cn)^{-\alpha} \qquad \forall \xi \in E^+(C).
\end{aligned} \tag{4.18}$$

The order of $N$ satisfies $\beta - 1 - 2\beta + 1 + \alpha = \alpha - \beta \geq 0$, which makes the right-hand side of (4.18) bounded from below by some $\varepsilon > 0$. Thus, condition (2.15) is fulfilled and Corollary 2.3 implies almost sure explosion on the set of all trajectories $(\zeta_0, \zeta_1, \dots)$ living in $E^+(C)$. Note that any $\xi \in \hat{E}$ satisfies $\xi \in E^+(C)$ for sufficiently large $C$, and

$$\zeta_0 \in E^+(C) \implies \zeta_k \in E^+(C) \qquad \forall k \text{ a.s.},$$

since mass is conserved and the number of particles increases. Thus, the assertion follows from Remark 4.1. $\square$

4.2. *Explosion of the mass flow process.* We consider the process with the kernel (3.29). In this model the growth of the size of individual particles is not bounded so that explosion is possible in the pure coagulation case. Theorem 4.4 gives a surprisingly easy and rather complete solution to this problem. As to pure fragmentation, again only the continuous case is of interest in the sense of explosion. Fragmentation in the mass flow process does not lead to a blow-up of the number of particles, but only to a decrease of their sizes. Theorem 4.7 provides sufficient conditions for explosion in the pure fragmentation case.

THEOREM 4.4. *Consider the mass flow kernel* (3.29) *in the pure coagulation case* ($S = 0, e = 0, F = 0$). *Assume*

$$(4.19) \qquad K(x,y) \geq \bar{K}(x,y),$$

*where $\bar{K}$ is homogeneous with exponent $\alpha > 1$, that is,*

$$(4.20) \qquad \bar{K}(cx, cy) = c^\alpha \bar{K}(x,y) \qquad \forall c > 0, x, y \in X,$$

*and such that $\bar{K}(1,1) > 0$. Then the mass flow process explodes almost surely, for any initial distribution on the set $\hat{E} = \{\xi \in E : N \geq 1\}$.*

PROOF. In order to make use of Corollary 2.3, we introduce the sets

$$E^+(C, L) = \left\{ \xi \in E : \min_{i=1,\dots,N} x_i \geq C, N = L \right\}, \qquad C > 0, L = 1, 2, \dots,$$

and the bounded measurable function

$$(4.21) \qquad \eta(\xi) = \begin{cases} \int_X H(x)\xi(dx), & \text{if } \xi \in E^+(C,L), \\ 0, & \text{otherwise}, \end{cases} \qquad \xi \in E,$$

where

$$H(x) = -x^{-\beta} \qquad \text{for some } 0 < \beta \leq \alpha - 1.$$

Consider trajectories $(\zeta_0, \zeta_1, \dots)$. Note that

$$(4.22) \qquad \zeta_0 \in E^+(C, L) \implies \zeta_k \in E^+(C, L) \qquad \forall k \text{ a.s.},$$

since the sizes of individual particles increase and the number of particles is conserved. According to (3.33), (3.14), (4.19) and (4.20), one obtains

$$\int_E [\eta(\xi_1) - \eta(\xi)]\tilde{q}(\xi, d\xi_1)$$

$$= \frac{1}{n^2} \sum_{i,j=1}^{L} [x_i^{-\beta} - (x_i + x_j)^{-\beta}] \frac{K_{\text{sym}}(x_i, x_j)}{x_j}$$

$$\geq \frac{1}{n^2} \sum_{i=1}^{L} [x_i^{-\beta} - (2x_i)^{-\beta}] x_i^{\alpha-1} \bar{K}(1,1)$$

$$= \frac{1 - 2^{-\beta}}{n^2} \bar{K}(1,1) \sum_{i=1}^{L} x_i^{\alpha-1-\beta}$$

$$\geq \frac{1 - 2^{-\beta}}{n^2} \bar{K}(1,1) L C^{\alpha-1-\beta} > 0 \qquad \forall \xi \in E^+(C, L).$$

Thus, condition (2.15) is fulfilled and Corollary 2.3 implies almost sure explosion on the set of all trajectories living in $E^+(C, L)$. Note that any $\xi \in \hat{E}$ satisfies $\xi \in E^+(C, L)$ for sufficiently small $C$ and some $L$. Taking into account (4.22), the assertion follows from Remark 4.1. □

EXAMPLE 4.5. Consider the mass flow coagulation process MF(1) starting with one particle of size one, that is, $\zeta_0 = \delta_1, n = 1$. In this case the Markov chain is deterministic,

$$\zeta_k = \delta_{x_k}, \qquad x_k = 2^k, \qquad k = 0, 1, \ldots,$$

and the sequence of waiting time parameters takes the form [cf. (3.32)]

$$\tilde{\lambda}(\zeta_k) = \frac{K(2^k, 2^k)}{2^k}.$$

Assuming $K = \bar{K}$ and using (4.20), one obtains $\tilde{\lambda}(\zeta_k) = 2^{k(\alpha-1)}, k = 0, 1, \ldots,$ so that

$$\sum_{k=0}^{\infty} \frac{1}{\tilde{\lambda}(\zeta_k)} < \infty \iff \alpha > 1.$$

EXAMPLE 4.6. Consider the mass flow coagulation process MF(2) starting with two particles of size one, that is, $\zeta_0 = \delta_1, n = 2$. In this case the Markov chain takes the form

$$\zeta_k = \tfrac{1}{2}[\delta_{x_k} + \delta_{y_k}], \qquad k = 0, 1, \ldots, x_0 = y_0 = 1,$$

and the sequence of waiting time parameters is [cf. (3.32)]

$$(4.23) \quad \tilde{\lambda}(\zeta_k) = \frac{1}{2}\left[\frac{K(x_k, x_k)}{x_k} + \frac{K(y_k, y_k)}{y_k} + \frac{K(x_k, y_k)}{y_k} + \frac{K(x_k, y_k)}{x_k}\right].$$

If

$$(4.24) \quad K(x,y) \geq (xy)^\beta \quad \forall x, y \in X, \text{ for some } \beta,$$

then one obtains

$$(4.25) \quad \tilde{\lambda}(\zeta_k) \geq \tfrac{1}{2}[x_k^{2\beta-1} + y_k^{2\beta-1} + x_k^\beta y_k^{\beta-1} + x_k^{\beta-1} y_k^\beta].$$

We consider several special trajectories, for which the explosion property can be checked explicitly.

First there is a trajectory of "fastest growth,"

$$(4.26) \quad x_k = 2^k, \quad y_k = 1, \quad k = 0, 1, \ldots,$$

when each jump consists in doubling the first particle. If the coagulation kernel satisfies (4.24), then one obtains from (4.25)

$$(4.27) \quad \tilde{\lambda}(\zeta_k) \geq \tfrac{1}{2} x_k^\beta$$

and

$$\sum_{k=0}^\infty \frac{1}{\tilde{\lambda}(\zeta_k)} \leq \sum_{k=0}^\infty \frac{2}{2^{\beta k}} < \infty \quad \text{if } \beta > 0.$$

Thus, there is explosion on the trajectory (4.26) even for many nongelling kernels.

Next we consider a trajectory of "second fastest growth,"

$$(4.28) \quad (1,1) \to (1,2) \to (3,2) \to (3,5) \to (8,5) \to (8,13) \to \cdots,$$

when self-interaction is avoided and alternatingly either the first particle is added to the second or vice versa. This trajectory is related to the Fibonacci numbers

$$\frac{(1+\sqrt{5})^k - (1-\sqrt{5})^k}{2^k \sqrt{5}} = \text{Integer}\left[\frac{1}{\sqrt{5}}\left(\frac{1+\sqrt{5}}{2}\right)^k\right], \quad k = 1, 2, \ldots,$$

where Integer[$a$] denotes the nearest integer of the number $a$. Here both particles $x_k, y_k$ grow as $C^k$ for some $C > 1$. If the coagulation kernel satisfies (4.24), then (4.25) implies $\tilde{\lambda}(\zeta_k) \geq C^{(2\beta-1)k}$ so that there is explosion on the trajectory (4.28) in the gelling case $\beta > 1/2$.

Finally, we consider a trajectory of "slowest growth,"

$$(4.29) \quad x_k = 1 + k, \quad y_k = 1, \quad k = 0, 1, \ldots,$$

when each jump consists in adding the second particle to the first one. If the coagulation kernel satisfies (4.24), then the estimate [cf. (4.27)]

$$\sum_{k=0}^{\infty} \frac{1}{\tilde{\lambda}(\zeta_k)} \leq \sum_{k=0}^{\infty} \frac{2}{(1+k)^\beta}$$

implies explosion in the case $\beta > 1$. If, however, $K(x,y) \leq xy$, $\forall x, y \in X$, then one obtains [cf. (4.23)] $\tilde{\lambda}(\zeta_k) \leq x_k + y_k = k+2$ and

$$\sum_{k=0}^{\infty} \frac{1}{\tilde{\lambda}(\zeta_k)} \geq \sum_{k=0}^{\infty} \frac{1}{k+2} = \infty.$$

Thus, there is no explosion on the trajectory (4.29) even for many gelling kernels.

THEOREM 4.7. *Consider the mass flow kernel* (3.29) *in the pure fragmentation case* $(S=0, e=0, K=0)$. *Assume*

$$(4.30) \qquad \tilde{F}(x, X) \geq \frac{C_F}{x^\alpha} \qquad \forall x \in X, \text{ for some } C_F > 0, \alpha > 0,$$

*and*

$$(4.31) \qquad \frac{1}{\tilde{F}(x,X)} \int_X \left(\frac{y}{x}\right)^\alpha \tilde{F}(x, dy) \leq \gamma < 1 \qquad \forall x \in X.$$

*Then the mass flow process explodes almost surely, for any initial distribution on the set* $\hat{E} = \{\xi \in E : N \geq 1\}$.

PROOF. In order to make use of Corollary 2.3, we consider the sets

$$E^+(C, L) = \left\{\xi \in E : \max_{i=1,\ldots,N} x_i \leq C, N = L\right\}, \qquad C > 0, L = 1, 2, \ldots,$$

and the bounded measurable function (4.21) with $H(x) = -x^\alpha$. Consider trajectories $(\zeta_0, \zeta_1, \ldots)$. Note that

$$(4.32) \qquad \zeta_0 \in E^+(C, L) \implies \zeta_k \in E^+(C, L) \qquad \forall k \text{ a.s.},$$

since the sizes of individual particles decrease and the number of particles is conserved. According to (3.33), (4.30) and (4.31), one obtains

$$\int_E [\eta(\xi_1) - \eta(\xi)] \tilde{q}(\xi, d\xi_1)$$

$$= \frac{1}{n} \sum_{i=1}^L \int_X [x_i^\alpha - y^\alpha] \tilde{F}(x_i, dy)$$

$$= \frac{1}{n} \sum_{i=1}^L x_i^\alpha \tilde{F}(x_i, X) \left[1 - \frac{1}{\tilde{F}(x_i, X)} \int_X \left(\frac{y}{x_i}\right)^\alpha \tilde{F}(x_i, dy)\right]$$

$$\geq \frac{C_F(1-\gamma)L}{n} > 0 \qquad \forall \xi \in E^+(C, L).$$

Thus, condition (2.15) is fulfilled and Corollary 2.3 implies almost sure explosion on the set of all trajectories living in $E^+(C,L)$. Note that any $\xi \in \hat{E}$ satisfies $\xi \in E^+(C,L)$ for sufficiently large $C$ and some $L$. Taking into account (4.32), the assertion follows from Remark 4.1. $\square$

In terms of the direct simulation fragmentation kernel $F$ [cf. (3.25) and (3.26)], assumption (4.30) takes the form (4.13), while assumption (4.31) takes the form

$$(4.33) \quad \frac{F^{(1)}_{\text{sym}}(x,X)}{F(x,Z)} \frac{1}{F^{(1)}_{\text{sym}}(x,X)} \int_X \left(\frac{y}{x}\right)^{\alpha+1} F^{(1)}_{\text{sym}}(x,dy) \leq \gamma < 1 \qquad \forall\, x \in X,$$

where the first factor represents the average number of fragments [cf. (3.21)]. The following examples illustrate assumption (4.33).

EXAMPLE 4.8. Consider the case of uniform binary fragmentation [cf. (3.23)]

$$F^{(1)}(x,dy) = F^{(1)}_{\text{sym}}(x,dy) = \bar{F}(x)\frac{1}{x}\mathbb{1}_{(0,x)}(y)\,dy, \qquad x \in X = (0,\infty),$$

where the function $\bar{F}$ determines the waiting time parameter. Assumption (4.33) is fulfilled, since

$$\frac{2}{x^{\alpha+2}}\int_0^x y^{\alpha+1}\,dy = \frac{2}{\alpha+2} < 1 \qquad \forall\, \alpha > 0.$$

EXAMPLE 4.9. Consider the case of deterministic binary fragmentation [cf. (3.23)]

$$(4.34) \qquad F^{(1)}(x,dy) = \bar{F}(x)\delta_{\kappa(x)}(dy),$$

where

$$\kappa(x) \in (0,x) \qquad \forall\, x \in X = (0,\infty),$$

and the function $\bar{F}$ determines the waiting time parameter. Note that [cf. (3.24)]

$$(4.35)\quad \begin{aligned} F^{(1)}_{\text{sym}}(x,dy) &= \frac{\bar{F}(x)}{2}\left[\delta_{\kappa(x)}(dy) + \int_X \delta_{\kappa(x)}(dx_1)\delta_{x-x_1}(dy)\right] \\ &= \frac{\bar{F}(x)}{2}[\delta_{\kappa(x)}(dy) + \delta_{x-\kappa(x)}(dy)]. \end{aligned}$$

Thus, assumption (4.33) takes the form

$$(4.36) \qquad \left(\frac{\kappa(x)}{x}\right)^{\alpha+1} + \left(\frac{x-\kappa(x)}{x}\right)^{\alpha+1} \leq \gamma < 1 \qquad \forall\, x \in X,$$

and is fulfilled, if

$$\gamma_1 x \leq \kappa(x) \leq \gamma_2 x \qquad \forall\, x \in X, \text{ for some } 0 < \gamma_1 \leq \gamma_2 < 1,$$

since

$$\left(\frac{\kappa(x)}{x}\right)^{\alpha+1} + \left(\frac{x-\kappa(x)}{x}\right)^{\alpha+1} \leq \frac{\kappa(x)}{x}\gamma_2^\alpha + \frac{x-\kappa(x)}{x}(1-\gamma_1)^\alpha$$

$$\leq \max(\gamma_2^\alpha, (1-\gamma_1)^\alpha) < 1.$$

In the last two examples we consider the case of deterministic binary fragmentation (4.34) and the mass flow fragmentation process with $n=1$ and one initial particle of size $x_0$. Note that [cf. (3.25) and (4.35)]

$$\tilde{F}(x, dy) = \frac{\bar{F}(x)}{2x}[y\delta_{\kappa(x)}(dy) + y\delta_{x-\kappa(x)}(dy)]$$

$$= \frac{\bar{F}(x)}{2x}[\kappa(x)\delta_{\kappa(x)}(dy) + [x-\kappa(x)]\delta_{x-\kappa(x)}(dy)]$$

and

$$\frac{1}{\tilde{F}(x,X)}\tilde{F}(x,dy) = \frac{\kappa(x)}{x}\delta_{\kappa(x)}(dy) + \left[1 - \frac{\kappa(x)}{x}\right]\delta_{x-\kappa(x)}(dy).$$

Thus, the Markov chain takes the form $\zeta_k = \delta_{x_k}, k = 0, 1, \ldots,$ where

$$(4.37) \qquad x_{k+1} = \begin{cases} \kappa(x_k), & \text{with probability } \dfrac{\kappa(x_k)}{x_k}, \\ x_k - \kappa(x_k), & \text{with probability } 1 - \dfrac{\kappa(x_k)}{x_k}. \end{cases}$$

The sequence of waiting time parameters is [cf. (3.32)]

$$(4.38) \qquad \tilde{\lambda}(\zeta_k) = \tilde{F}(x_k, X) = \tfrac{1}{2}\bar{F}(x_k).$$

EXAMPLE 4.10. Consider the case (4.34) with the function

$$(4.39) \qquad \kappa(x) = \begin{cases} \dfrac{x}{2} + \dfrac{1}{4}, & \text{if } x > \tfrac{1}{2}, \\ \dfrac{x}{2}, & \text{otherwise.} \end{cases}$$

Note that assumption (4.33) [cf. (4.36)] is violated, since $\lim_{\varepsilon \to 0}\kappa(\tfrac{1}{2}+\varepsilon) = \tfrac{1}{2}$. If $x_0 > \tfrac{1}{2}$, then the sequence

$$\eta_k = \frac{2x_0 + 2^k - 1}{2^{k+1}}, \qquad k = 0, 1, \ldots,$$

satisfies

$$\kappa(\eta_k) = \frac{2x_0 + 2^k - 1}{2^{k+2}} + \frac{1}{4} = \eta_{k+1}$$

and corresponds to the trajectory of "slowest decrease" [cf. (4.37)]. There is no explosion on this trajectory, since it does not reach zero. However, it has a nonzero probability, namely,

$$(4.40) \qquad \lim_{k\to\infty} \frac{\kappa(\eta_0)}{\eta_0}\frac{\kappa(\eta_1)}{\eta_1}\cdots\frac{\kappa(\eta_k)}{\eta_k} = \frac{1}{x_0}\lim_{k\to\infty} \kappa(\eta_k) = \frac{1}{2x_0}.$$

This example illustrates that some additional restriction like assumption (4.31) in Theorem 4.7 cannot be avoided.

EXAMPLE 4.11. Consider the case (4.34) with the function $\kappa(x) = \frac{x}{2}$ and $x_0 = 1$. According to (4.37), the trajectory of the Markov chain is deterministic, with

$$x_k = 2^{-k}, \qquad k = 0, 1, \ldots.$$

For the choice $\bar{F}(x) = -\log x + 1$, one obtains [cf. (4.38)]

$$\sum_{k=0}^{\infty} \frac{1}{\tilde{\lambda}(\zeta_k)} = \sum_{k=0}^{\infty} \frac{2}{\bar{F}(x_k)} = \frac{2}{\log 2} \sum_{k=1}^{\infty} \frac{1}{k} = \infty$$

so that there is no explosion. This example illustrates that assumption (4.30) in Theorem 4.7 cannot be replaced by some arbitrarily slow growth at zero.

4.3. *Comments.* In the pure coagulation case there is no explosion in the direct simulation model. According to Theorem 4.4, there is explosion in the mass flow model for a rather wide class of gelling coagulation kernels.

In the pure fragmentation case there is explosion in both models. Theorems 4.3 and 4.7 cover wide classes of unbounded (at zero) fragmentation kernels. The sufficient conditions for the mass flow model are stronger than those for the direct simulation model. Example 4.10 shows that there are fragmentation kernels for which the direct simulation model explodes almost surely, while the mass flow model does not. Indeed, choose the kernel (4.34) with $\bar{F}(x) = x^{-\alpha}, \alpha > 0$, and $\kappa$ given in (4.39). Then the direct simulation process explodes almost surely for any initial value $x_0$, according to Theorem 4.3. The mass flow process explodes with probability $1 - \frac{1}{2x_0}$ for any $x_0 > \frac{1}{2}$, according to (4.40).

In the direct simulation fragmentation model explosion implies that the number of particles in the system reaches infinity in finite time (each jump increases this number). Due to mass conservation, this means that (at the explosion time) infinitely many particles are below any given size $\varepsilon$, creating "dust." In the mass flow fragmentation model explosion occurs due to a fast approach to zero of a single particle. This explains what is going on in the example described above. In the direct simulation case each jump creates two fragments and some of them reach the interval $(0, \frac{1}{2}]$ leading to explosion.

Fragments staying bigger than $\frac{1}{2}$ do not avoid the fact of explosion. In the mass flow model there is only one particle, which stays above $\frac{1}{2}$ with positive probability.

Recently the interest in fragmentation processes has considerably increased (cf., e.g., [4, 15]). In this context the occurrence of small particles (dust) has been studied using different models.

**5. Concluding remarks.** In this paper explosion criteria for jump processes with an arbitrary locally compact separable metric state space were established. These results are of independent interest. As an illustration, the general criteria were applied to stochastic coagulation–fragmentation models. The corresponding results (Theorems 4.3, 4.4 and 4.7) cover a wide range of coagulation and fragmentation kernels, for which phase transitions to infinitely large (gel) or infinitely small (dust) particles are known in the context of macroscopic equations. The proofs of these theorems are rather short, which illustrates the efficiency of the general explosion criteria. In particular, the results of Theorems 4.3 and 4.4 would hardly be available on the basis of the previously known criteria for one-dimensional processes. It might be of interest to extend the explosion results, which were obtained either for pure coagulation or for pure fragmentation, to more general situations combining both processes and including source and efflux terms. For this purpose, perhaps more sophisticated choices of the test function $\eta$ in the criteria will be required.

A challenging problem for future research is the study of the limiting behavior ($n \to \infty$) of the stochastic coagulation–fragmentation models with explosion. In particular, the second (more difficult) part of the conjecture mentioned in Section 1 is still open. Convergence of a truncated mass flow model (particles exceeding a certain level are removed from the system) to the solution of the Smoluchowski equation has been studied in [9]. In general, the continuation of jump processes after the explosion point is a delicate problem. In the context of the coagulation–fragmentation models based on particles, there is a natural way of continuation—particles reaching 0 or $\infty$ are simply removed from the system, while the others continue their evolution according to the previous rules.

The conservation property of the mass flow equation [cf. (3.27) without source and efflux terms and for a constant test function] allows one to interpret its solution as a probability measure and to construct a related nonlinear Markov process. Such a process is determined by some stochastic equation with coefficients depending on the law of the solution. This approach has been carried out in [7] (pure coagulation), [17] (including discrete fragmentation) and [13] (including continuous fragmentation). It would be of interest to clarify how the explosion phenomena observed in finite particle systems are represented in these limiting processes.

It is remarkable that some explosion phenomenon is probably recovered also in the direct simulation coagulation model in the limit $n \to \infty$. A result concerning explosion of an appropriately scaled tagged particle in the discrete case $K(x,y) = (xy)^\alpha$ for $1/2 < \alpha \leq 1$ was announced in [21], but unfortunately has not been published so far.

**Acknowledgments.** Most of the results of this paper were obtained during my stay at the Isaac Newton Institute for Mathematical Sciences (Cambridge, UK) in October–December 2003. I am very grateful to the organizers of the scientific programme "Interaction and Growth in Complex Stochastic Systems" and to the staff of the Newton Institute for creating a very encouraging and pleasant research environment. In particular, I would like to thank James Norris for the invitation and several helpful discussions.

WEIERSTRASS INSTITUTE FOR APPLIED
ANALYSIS AND STOCHASTICS
MOHRENSTRAßE 39
D-10117 BERLIN
GERMANY
E-MAIL: wagner@wias-berlin.de